\documentstyle[amssymb,amsfonts]{amsart}

\def\Z{{\hbox{\bf Z}}}

\def\eps{\varepsilon}

\newenvironment{proof}{\noindent {\bf Proof} }{\endprf\par}
\def \endprf{\hfill  {\vrule height6pt width6pt depth0pt}\medskip}
\def\emph#1{{\it #1}}
\def\textbf#1{{\bf #1}}

\theoremstyle{plain}
  \newtheorem{theorem}[subsection]{Theorem}

  \newtheorem{lemma}[subsection]{Lemma}
  \newtheorem{corollary}[subsection]{Corollary}

\theoremstyle{remark}
  \newtheorem{remark}[subsection]{Remark}

\theoremstyle{definition}

\include{psfig}

\begin{document}

\title{A remark on primality testing and decimal expansions}
\author{Terence Tao}
\address{Department of Mathematics, UCLA, Los Angeles CA 90095-1555}
\email{ tao@@math.ucla.edu}

\vspace{-0.3in}
\begin{abstract}
We show that for any fixed base $a$, a positive proportion of primes have the property that they become composite after altering any one of their digits in the base $a$ expansion; the case $a=2$ was already established by Cohen-Selfridge and Sun, using some covering congruence ideas of Erd\H{o}s.  Our method is slightly different, using a partially covering set of congruences followed by an application of the Selberg sieve upper bound.  As a consequence, it is not always possible to test whether a number is prime from its base $a$ expansion without reading all of its digits.  We also present some slight generalisations of these results.
\end{abstract}

\maketitle

\section{Introduction}

In 1950, Erd\H{o}s\cite{erdos} used the method of covering congruences to show that there exists an infinite arithmetic progression of odd integers $m$ with the property that $|m-2^i|$ is composite for every $i$.  Modifying this method, Cohen and Selfridge \cite{cs} exhibited an arithmetic progression of odd integers $m$ such that $|m-2^i|$ and $m+2^i$ are both composite for every $i$.  In \cite{sun}, Sun gave the explicit arithmetic progression $\{ m: m = M \mod \prod_{p \in {\mathcal P}} p\}$ with this property, where
$$ M := 47867742232066880047611079$$
and ${\mathcal P}$ is the finite set of primes
$$ {\mathcal P} := \{2,3,5,7,11,13,17,19,31,37,41,61,73,97,109,151,241,257,231\},$$
and noted that integers in this progression are in fact not of the form $\pm p^a \pm q^b$ for any primes $p,q$ and positive integers $a,b$.  Since $M$ is coprime to $\prod_{p \in {\mathcal P}} p$, we can apply the prime number theorem in arithmetic progressions (see e.g. \cite[Corollary 11.17]{mv}) to obtain the following immediate corollary:

\begin{corollary}\label{main}\cite{cs}, \cite{sun}  For all sufficiently large integers $n$, there exist at least $c 2^n/n$  primes $p$ between $2^{n-1}$ and $2^n$ such that the integers $p-2^i$ and $p+2^i$ are composite for every $0 \leq i < n-1$, where $c > 0$ is an absolute constant.
\end{corollary}

We remark that primes $p$ of the above form are initially rather rare; the first few primes of this form are 
$$1973, 3181, 3967, 4889, 8363, 8923, 11437, 12517, 14489,\ldots.$$
On the other hand, from Corollary \ref{main} and the prime number theorem we see that a positive proportion of the primes in fact lie on this sequence.

As an immediate corollary of Corollary \ref{main}, we see that for sufficiently large $n$, there exist $n$-bit integers $p$ which are prime, but such that any number formed from $p$ by switching one of the bits is not prime; the first few primes of this form are $127, 173, 191, 223, 233, \ldots$ (a slight variant of sequence A065092 in \cite{OEIS}, which is the subsequence in which $p+2^{n+1}$ is also required to be composite).  In other words, if we let $P_n: \{0,1\}^n \to \{0,1\}$ be the boolean function which returns $1$ if and only if the $n$-bit integer corresponding to the input $\{0,1\}^n$ is prime, then the \emph{sensitivity} $s(P_n)$ of $P_n$ is equal to $n$ for sufficiently large $n$.  Recall that the sensitivity (or \emph{critical complexity}) $s(B)$ of a Boolean function $B: \{0,1\}^n \to \{0,1\}$ is the largest integer $s$ for which there exists an input $x \in \{0,1\}^n$ such that $B(x) \neq B(x')$ for at least $s$ inputs $x'$ which are formed from $x$ by switching exactly one bit.  We remark that the lower bound $s(P_n) \geq \frac{1}{4}n + O(1)$ was previously established in \cite[p. 307]{crypt}.

If $p$ is as above, then clearly it is not possible for an algorithm to determine with absolute certainty whether $p$ is prime or not without inspecting all of the digits in the binary expansion.  In particular, any deterministic primality tester can require computational time at least logarithmic in the size of the number being tested, if that number is represented in binary.  For comparison, it was shown in \cite[Theorem 6]{mosc} that any recursive algorithm which can decide the primality of an $n$-bit integer using the operations $=$, $<$, $+$, $-$, $2 \cdot$, $\frac{1}{2} \cdot$, and parity, has time complexity at least $\frac{1}{4} n$.  We remark that for bounded depth circuits, much stronger lower bounds (of exponential type in $n$) on the spatial complexity are known; see \cite{ass}, \cite{woods}.

In this note we establish a similar result for general bases.  More precisely, we establish

\begin{theorem}\label{main2} Let $K \geq 1$ be an integer.  Then for all sufficiently large $N$, the number of primes $p$ between $N$ and $(1+\frac{1}{K})N$ such that $|k p \pm j a^i|$ is composite for all integers $1 \leq a, j, k \leq K$ and $1 \leq i \leq K \log N$ is at least $c_K \frac{N}{\log N}$ for some constant $c_K > 0$ depending only on $K$.
\end{theorem}

From this theorem we see that the above results for binary expansions are also valid in other bases as well.  For instance, applying this theorem with $K=10$ we conclude that a positive proportion of the primes have the property that if one changes any one of the digits in the base $10$ expansion, one necessarily obtains a composite number, and so any deterministic primality tester receiving the digits of this number as input must read all of these digits in order to determine its primality.  The first few such primes are $294001, 505447, 584141, \ldots$ (sequence A050249 from \cite{OEIS}).  The infinitude of this sequence was established previously by Erd\H{o}s \cite{orno}.

Our argument does not use a fully covering set of congruences.  Instead, one uses congruences modulo primes arising from Mersenne-type numbers (in which bases such as $a$ have an unexpectedly low order) to sieve out most of the quadruples $(a,j,k,i)$ appearing in the above theorem, leaving behind a small number which can be handled via standard upper bound sieves.  It seems to be difficult to establish this result without such a preliminary sieving step, since without such a sieving one would expect each $|kp \pm j a^i|$ to be prime with probability comparable to $\frac{1}{\log N}$, which makes it moderately unlikely (especially for large $K$) that the $|kp \pm j a^i|$ are composite for all $1 \leq a,j,k \leq K$ and $1 \leq i \leq K \log N$ for any given prime $p$.

The author is supported by NSF grant CCF-0649473 and a grant from the MacArthur Foundation.  The author is indebted to Yiannis Moschovakis for suggesting this question, and to Jens Kruse Andersen, Yong-Gao Chen, Bjorn Poonen, Florian Luca, Paul Pollack, Igor Shparlinski, Zhi-Wei Sun, and several anonymous commenters on my blog for helpful comments and references.

\section{Proof of Theorem \ref{main2}}

We now prove Theorem \ref{main2}.  Fix $K$.  We will need a large integer $M = M(K) \geq K$ to be chosen later.  We will then use this integer $M$ to generate a finite set ${\mathcal P}$ of primes, as follows:

\begin{lemma}  For any $M, K \geq 1$, there exists a finite set ${\mathcal P}$ of primes which can be partitioned into disjoint sets ${\mathcal P} = \bigcup_{2 \leq a \leq K} {\mathcal P}_a$, with the following properties:
\begin{itemize}
\item If $p \in {\mathcal P}_a$ for some $2 \leq a \leq K$, then there exists a prime $q_p$ such that
\begin{equation}\label{qp-large}
q_p \geq Mp
\end{equation}
and
\begin{equation}\label{ap-small}
a^p = 1 \hbox{ mod } q_p.
\end{equation}  
Furthermore, the primes $q_p$ for $p \in {\mathcal P}$ are all distinct.
\item For each $2 \leq a \leq K$, we have
\begin{equation}\label{pam}
\sum_{p \in {\mathcal P}_a} \frac{1}{p} \geq M.
\end{equation}
\end{itemize}
\end{lemma}

\begin{proof}
The claim is trivial for $K=1$, so assume inductively that $K \geq 2$ and that the claim has already been proven for $K-1$.  Thus we already have disjoint finite sets of primes ${\mathcal P}_1,\ldots,{\mathcal P}_{K-1}$ with the stated properties.

Let $W$ denote the product of all the numbers less than $K$ which are coprime to $K$, and let $A$ denote the multiplicative order of $K \hbox{ mod } W^K$.  Observe that if $p$ is a prime with $p = 1 \hbox{ mod } A$, then $K^p-1 = K-1 \hbox{ mod } W^K$.  In particular, if $q$ is any prime less than $K$, then $q$ can divide $K^p-1$ at most $K$ times (since $K-1$ is not a multiple of $q^K$, being a smaller integer).  As a consequence, we see that if $p$ is larger than some sufficiently large constant $C_K$, then the largest prime factor of $K^p-1$ is greater than $K$.

By the prime number theorem in arithmetic progressions (see e.g. \cite[Corollary 11.17]{mv}), the sum of reciprocals of primes equal to $p=1 \hbox{ mod } A$ is divergent.  From this and Corollary \ref{sophie} we may find a infinite collection of primes ${\mathcal P}'$ of primes $p=1 \hbox{ mod } A$ which are larger than $C_K$, disjoint from the finite sets ${\mathcal P}_1,\ldots,{\mathcal P}_{K-1}$, such that $\sum_{p \in {\mathcal P}'} \frac{1}{p} = \infty$, and such that $mp+1$ is composite for every $1 \leq m \leq M$.  For any $p$ in ${\mathcal P}'$, we set $q_p$ to be the largest prime factor of $K^p-1$.  Since $p > C_K$, we have $q_p > K$.  In particular, the multiplicative order of $K \hbox{ mod } q_p$ is exactly $p$, which forces all the $q_p$ to be distinct.
In particular, we can find a finite subset ${\mathcal P}_K$ of ${\mathcal P}'$ with $\sum_{p \in {\mathcal P}_K} \frac{1}{p} \geq M$ such that the values of $q_p$ for $p \in {\mathcal P}_K$ are distinct from all the values of $q_p$ already assigned to $p$ in ${\mathcal P}_1,\ldots,{\mathcal P}_{K-1}$.

From Fermat's little theorem we see that $p$ divides $q_p-1$ for all $p \in {\mathcal P}_K$.  On the other hand, we have $mp+1$ composite for every $1 \leq m \leq M$.  Thus $q_p \geq Mp$ as required.  Thus ${\mathcal P} := {\mathcal P}_1 \cup \ldots \cup {\mathcal P}_K$ obeys all the desired properties.
\end{proof}

\begin{remark} In \cite{shorey} it is shown that the largest prime factor of $2^p-1$ is at least $c p \log p$ for some absolute constant $c > 0$ (see also \cite{murata} for additional refinements and further discussion).  Slightly weaker results for more general bases can be found in \cite{lusz}.  By using these results one can avoid the use of Corollary \ref{sophie}.
\end{remark}

Henceforth we let ${\mathcal P} = {\mathcal P}_2 \cup \ldots \cup {\mathcal P}_K$, as well as the primes $q_p$ for $p \in P$ be as in the above lemma.   

We let $N$ be a sufficiently large integer parameter.
We use the asymptotic notation $o(1)$ to denote any quantity that goes to zero as $N \to \infty$ (with $K$, $M$, and ${\mathcal P}$ fixed), and similarly $X \ll Y$ or $X = O(Y)$ to denote the estimate $X \leq CY$ for some $C$ depending on $K$ but independent of $N$, $M$, ${\mathcal P}$.  We also write $X \sim Y$ for $X \ll Y \ll X$.

By reducing the sets ${\mathcal P}_a$ if necessary, we may assume from \eqref{pam} that
\begin{equation}\label{pam2}
\sum_{p \in {\mathcal P}_a} \frac{1}{p} \sim M.
\end{equation}
Let $S$ denote the finite set of pairs
$$ S := \{ (j,k) \in \Z^2: -K \leq j \leq K; 1 \leq k \leq K; j \neq 0 \}.$$
By \eqref{pam} and a simple greedy argument, we may partition ${\mathcal P}_a = \bigcup_{(j,k) \in S} {\mathcal P}_{a,j,k}$ in such a way that
\begin{equation}\label{prime-half}
\sum_{p \in {\mathcal P}_{a,j,k}} \frac{1}{p} \sim M
\end{equation}
for all $2 \leq a \leq K$ and $(j,k) \in S$.

Let $W$ be the quantity
$$W := \prod_{p \in {\mathcal P}} q_p.$$
By the Chinese remainder theorem, we can find $b$ coprime to $W$ such that $kb+j=0 \mod q_p$ for $p \in {\mathcal P}_{a,j,k}$, $2 \leq a \leq K$, and $(j,k) \in S$.  (Note from \eqref{qp-large} and the hypothesis $M \geq K$ that all integers between $1$ and $K$ are coprime to $W$.)

To establish Theorem \ref{main2}, it will suffice to show that the quantity
\begin{equation}\label{nnk}
\begin{split}
 \# \{ &N \leq m \leq (1+\frac{1}{K})N: m = b \hbox{ mod } W; m \hbox{ prime, but } \\
 &\quad |km + j a^i| \hbox{ composite for all } 0 \leq i < K \log N, 1 \leq a \leq K, (j,k) \in S \}
\end{split} 
\end{equation}
is $\gg N / \log N$.  Note that when $a=1$, the value of $i$ is irrelevant (and so can be set for instance to zero).  We can thus crudely bound \eqref{nnk} from below by
\begin{equation}\label{nnk2}
\eqref{nnk} \geq Q_N -  \sum_{a=2}^K \sum_{(j,k) \in S} \sum_{0 \leq i < K \log N} Q_{N,i,a,j,k} - 
\sum_{(j,k) \in S} Q_{N,0,1,j,k} - O(\log N)
\end{equation}
where
$$ Q_N := \# \{ N \leq m \leq (1+\frac{1}{K})N: m = b \hbox{ mod } W \}$$
and
$$ Q_{N,i,a,j,k} := \# \{N \leq m \leq (1+\frac{1}{K})N: m = b \hbox{ mod } W; m, |km \pm ja^i| \hbox{ both prime } \}.$$
(The $O(\log N)$ error arises from the small number of cases in which $|km+j a^i|$ is equal to zero or one.)

From the prime number theorem in arithmetic progressions (see e.g. \cite[Corollary 11.17]{mv}) we have
$$ Q_N \gg \frac{N}{\phi(W) \log N}$$
where
$$ \phi(W) = W \prod_{p \in {\mathcal P}} (1 - \frac{1}{q_p})$$
is the Euler totient function of $W$.  (More precise asymptotics for $Q_N$ are available, but we will not need them here.)  

From Corollary \ref{sieve} we have
\begin{equation}\label{qni}
 Q_{N,i,a,j,k} \ll \frac{N}{W \log^2 N} \prod_{p \in {\mathcal P}} (1 - \frac{1}{q_p})^{-2} 
\end{equation}
for all $1 \leq a \leq K$, $0 \leq i \leq K \log N$, and $(j,k) \in S$.  Applying this to dispose of the $Q_{N,0,1,j,k}$ terms in \eqref{nnk2}, we thus conclude that
\begin{equation}\label{nnk3}
\eqref{nnk} \gg
\frac{N}{W\log N} \prod_{p \in {\mathcal P}} (1 - \frac{1}{q_p})^{-1} - O( \sum_{a=2}^K \sum_{(j,k) \in S} \sum_{0 \leq i < K \log N} Q_{N,i,a,j,k} )
\end{equation}
when $N$ is sufficiently large.

Now suppose that $2 \leq a \leq K$ and $(j,k) \in S$.  Observe that if $i = 0 \mod p$ for any $p \in {\mathcal P}_{a,j,k}$, then $|km + ja^i|$ is divisible by $q_p$, and thus will prime for at most one value of $m$.  Thus (paying a negligible factor of $O(\log N)$) we may restrict attention to those $0 \leq i < n-1$ such that $i \neq 0 \mod p$ for every $p \in {\mathcal P}_{a,j,k}$.  By the Chinese remainder theorem, we see that the number of such $i$ is $O( \log N \prod_{p \in {\mathcal P}_{a,j,k}} (1 - \frac{1}{p}) )$.  Using the approximations
$$ \prod_{n \in A} (1 - \frac{1}{n}) \sim \exp( - \sum_{n \in A} \frac{1}{n} )$$
which are valid for any finite set $A$ (since $\sum_{n\in A} \frac{1}{n^2} = O(1)$), we conclude from the above discussion and \eqref{qni} that
$$ \sum_{0 \leq i < K \log N} Q_{N,i,a,j,k} \ll \frac{N}{W \log^2 N} 
\left(\prod_{p \in {\mathcal P}} (1 - \frac{1}{q_p})^{-1}\right)
\exp( \sum_{p \in {\mathcal P}} \frac{1}{q_p} - \sum_{p \in {\mathcal P}_{a,j,k}} \frac{1}{p} ).$$
But from \eqref{pam}, \eqref{qp-large} we have $\sum_{p \in {\mathcal P}} \frac{1}{q_p} = O(1)$, while from \eqref{prime-half} we have $\sum_{p \in {\mathcal P}_{a,j,k}} \frac{1}{p} \gg M$.  Inserting all these bounds into \eqref{nnk3}, we conclude
$$ \eqref{nnk} \gg \frac{N}{W\log N} ( 1 - O( \exp( - c M) ) ) \left(\prod_{p \in {\mathcal P}} (1 - \frac{1}{q_p})^{-1}\right) $$
where $c > 0$ depends on $K$ but not on $M$.  Taking $M$ sufficiently large depending on $K$, we obtain the claim.

\section{Remarks}

An inspection of the proof of Theorem \ref{main2} allows one to establish a strengthened version in which the numbers $|k p \pm j a^i|$ are not only composite, but they also contain at least two distinct prime factors greater than $K$.  More precisely, the cases in which $|kp \pm j a^i|$ is the product of a prime power $q^b$ and some primes less than or equal to $K$ can be disposed of by suitable variants of Corollary \ref{sieve} (and in the case $b \geq 2$, the total contribution here is $O(\sqrt{N})$ which is easily discarded); we omit the details.  Recently in \cite{pan2}, it was shown that one can in fact ensure that the numbers $kp \pm ja^i$ contain $\gg (\log \log N)^{1/3-\eps}$ prime factors each for any fixed $\eps$.

In a somewhat different direction, it should also be possible to strengthen the conclusion of Theorem \ref{main2} to assert that $|kp \pm ja^i + l|$ is composite for all $l$ in some set $L = L_N \subset \{-KN,\ldots,KN\}$ of cardinality at most $K$.  A new difficulty arises here due to an additional factor of $\prod_{p | \pm ja^i + l; p /\!\!| W} (1-\frac{1}{p})^{-1}$ arising from the use of Corollary \ref{sieve}, but it seems likely that this quantity should be bounded for the overwhelming majority of values of $a,i,j,l$, which should allow one to continue the argument; we will not pursue this matter here.  If one is able to carry out this generalisation, one should be able to obtain the conclusion that for any base $a \geq 2$ and any $r \geq 1$, a positive proportion of the primes $p$ have the property that if one modifies any single one of its digits in the base $a$ expansion, \emph{and} appends or deletes up to $r$ digits to the end and/or beginning of the digit string, one necessarily obtains a composite number.  

In a similar spirit, it was recently established in \cite{kozek} that there exist infinitely many composite numbers coprime to which remain composite after inserting a single digit in their base $10$ expansion.  It seems likely that one should now also be able to find infinitely many prime numbers with the same property (i.e. they become composite after inserting any digit at any place).

In all of the above results, the total number of possible modifications of the digit string remains comparable to $\log p$ and so the cases in which a number is unexpectedly prime can be handled by the upper bound sieve after performing the preliminary sieving to eliminate most of the cases.  The problem becomes significantly more difficult, however, if one asks that the number $p$ become composite after allowing one to modify any \emph{two} of the digits in the digit string, as the number of possible modifications is now comparable to $\log^2 p$. Indeed, standard heuristics from the prime tuples conjecture \cite{hardy-littlewood} now lead one to predict that for a sufficiently large base, there should only be finitely many numbers of this form, although there is a slim chance (especially in small bases) that Mersenne-type primes provide enough congruences to fully cover all the modifications for primes in a certain infinite arithmetic progression, as was the case with Theorem \ref{main}.  We remark that in \cite{yuan} it was shown that there are infinitely many integers $n$ such that $n - 2^a - 2^b$ is not a prime power for any $a, b$ (an earlier result in \cite{crocker} establishes the weaker statement with ``prime power'' replaced by ``prime'').	The base $2$ was generalised to other bases recently in \cite{chen}, and lower bounds on the density of such integers was obtained in \cite{chen} and \cite{pan} (the latter result using the methods in this paper).

Using the circle method and bounds on prime exponential sums, there are several further results known relating primes to binary digits, or to powers of $2$.  For instance, in \cite{harman} the distribution of a bounded number of fixed digits of a large prime was studied.  In \cite{mr} it was shown that the binary digit sum of a large prime was equally likely to be even as it was to be odd.  In a slightly different direction, it was shown in \cite{hp} that all sufficiently large even numbers are the sum of two primes, together with at most $13$ powers of two.

\appendix

\section{Some sieve theory}

We recall the following standard application of the Selberg sieve to twin prime type problems:

\begin{theorem}[Selberg sieve upper bound]\label{Selberg}  Suppose that $y \geq 4$, and let $P := \prod_{p < \sqrt{y}} p$.  Let ${\mathcal B}(p)$ be the union of $b(p)$ arithmetic progressions with common difference $p$, and put ${\mathcal B} := \bigcup_{p|P} {\mathcal B}(p)$.  If $b(2) \leq 1$ and $b(p) \leq 2$ for $p > 2$, then the number of integers $0 \leq r \leq y$ such that $r \not \in {\mathcal B}$ is 
$$ \ll \frac{y}{\log^2 y} \prod_{p|P} (1-\frac{b(p)}{p}) (1 - \frac{1}{p})^{-2}.$$
\end{theorem}

\begin{proof}  See \cite[Theorem 3.13]{mv}.  As shown in that reference, one can in fact replace the implied constant with $8 + O( \frac{\log \log y}{\log y})$, but we will not need this improvement here.
\end{proof}

\begin{corollary}\label{sieve}  Let $x, W, b, \geq 1$ be integers with $W$ even, and let $h, k$ be non-zero integers.  Then if $x$ is sufficiently large depending on $W, b$, we have
\begin{align*}
\# \{ 0 < m \leq x: m = b &\hbox{ mod } W; m, |km+h| \hbox{ both } prime \}\\
&\ll \frac{x}{W \log^2 x} \left(\prod_{p|W} (1 - \frac{1}{p})^{-2}\right) \left(\prod_{p|h; p /\!\!| W} (1 - \frac{1}{p})^{-1}\right)
\end{align*}
where the implied constant can depend on $k$.
\end{corollary}

\begin{proof}  By reversing the signs of $k$ and $h$ if necessary, and increasing the size of the implied constant by a factor of $2$ if necessary, we may replace $|km+h|$ by $km+h$. We may assume that $b$ and $kb+h$ are both coprime to $W$, otherwise the number of $m$ for which $km, km+h$ are both prime is bounded uniformly in $x$ and the claim is trivial.  For similar reasons we may assume that $k$ and $h$ are coprime.  Write $m = Wr+b$ and $y := x/W$, thus $0 \leq r \leq y$.  We can restrict attention to the case $r > \sqrt{y}$, since the case $r \leq \sqrt{y}$ only contributes $O(\sqrt{y})$ elements which is acceptable.  If $p \leq \sqrt{y}$ is a prime, then the constraints that $m$ and $m+h$ both be prime force $Wr+b$ and $kWr+kb+h$ to both be coprime to $p$.  If $p|W$, then this condition is vacuous; if $p|h$, $p /\!\!| W$, and $p /\!\!| k$, then this excludes one residue class modulo $p$ from the space of possible $r$'s; and if $p /\!\!| h$, $p /\!\!| W$ and $p /\!\!| k$ then this excludes two residue classes modulo $p$ from the space of possible $r$'s.  Finally, if $p/\!\!|W$ and $p|k$ then either one or two residue classes modulo $p$ are excluded.  The claim now follows from Theorem \ref{Selberg} (note that $\log x$ is comparable to $\log y$ for $x$ large enough, and that $\prod_p (1-\frac{2}{p}) (1-\frac{1}{p})^{-2}$ is comparable to $1$).
\end{proof}

\begin{corollary}[Brun's theorem]\label{sophie}  Let $m, j$ be any positive integers.  Then the sum of reciprocals of the primes $p$ for which $mp+j$ is also prime is convergent.
\end{corollary}

\begin{proof}  By Corollary \ref{sieve}, the number of primes of the above form which are less than $x$ is $O( \frac{x}{\log^2 x})$ (where the implied constant can depend on $m$).  The claim easily follows.
\end{proof}

\end{document}